\documentclass{elegantpaper}

\usepackage{amsmath}
\usepackage{amsfonts}
\usepackage{amssymb}
\usepackage{mathrsfs}
\usepackage{graphicx}
\usepackage{xcolor}

\newcommand  \bglb {\big (}
\newcommand  \bgrb {\big )}

\newcommand  \bgls {\big [}
\newcommand  \bgrs {\big ]}

\newcommand   \Integers {\mathbb Z}

\newcommand   \primeInt {\mathtt{q}}
\newcommand   \fieldChar {\mathtt{\,p\,}}
\newcommand   \finiteIntegerField {\basel{\mathbb Z}{\fieldChar}\,} 
\newcommand   \scalars  {\mathbb F}

\newcommand   \extensionField  {\mathbb {\widehat F}}

\renewcommand   \gcd  {~\textsf{gcd}~ }

\newcommand  \singlevariablepolynomials[2] {#1{\mathbf{[}#2\mathbf{]}}}
\newcommand  \basel[2]{#1_{_{#2}}}

\newcommand  \lspace  {\hspace*{-1.5mm}}

\newcommand  \tab  {\hspace*{0.5cm}}
\newcommand  \ltab  {\hspace*{-0.5cm}}

\newcommand  \shiftright  {\hspace*{2.0cm}}
\newcommand  \shiftleft  {\hspace*{-2.0cm}}

\definecolor{titlecolor}{RGB}{144,48,48}
\definecolor{authorname}{RGB}{16,96,16}%
\definecolor{addrscolor}{RGB}{60,113,183}
\definecolor{secheader}{RGB}{16,32,128}%
\definecolor{refscolor}{RGB}{16,64,128}%

\begin{document}

\title{\textcolor{titlecolor}{\bf{Optimal Normal Bases Over Finite Fields}}}

\author{\\
\\
{\textcolor{authorname}{\bf{Duggirala Meher Krishna}}}\\
{\textcolor{addrscolor}{\small{\bf{Gayatri Vidya Parishad College of Engineering (Autonomous)}}}} \\
{\textcolor{addrscolor}{\small{\bf{Madhurawada, VISAKHAPATNAM -- 530 048, Andhra Pradesh, India}}}} \\
{\textcolor{addrscolor}{\small{\bf{E-mail ~: \tab duggiralameherkrishna@gmail.com}}}}\\
 \\
 {\textcolor{addrscolor}{and}} \\
 \\
{\textcolor{authorname}{\bf{Duggirala Ravi}}}\\
{\textcolor{addrscolor}{\small{\bf{Gayatri Vidya Parishad College of Engineering (Autonomous)}}}} \\
{\textcolor{addrscolor}{\small{\bf{Madhurawada, VISAKHAPATNAM -- 530 048, Andhra Pradesh, India}}}} \\
\shiftleft \tab {\textcolor{addrscolor}{\small{\bf{E-mail ~: \tab ravi@gvpce.ac.in; \tab duggirala.ravi@yahoo.com};}}} \\
{\textcolor{addrscolor}{\small{\shiftright \ltab \ltab \bf{ duggirala.ravi@rediffmail.com; \tab drdravi2000@yahoo.com}}}}
\\
\\
}

\date{}

\maketitle

\textcolor{secheader}{
\begin{abstract}
% Text of abstract
In this paper, a method for constructing a near optimal normal basis for algebraic extensions of a finite field is described. In each extension, except for the squares of basis elements, the product of two distinct normal basis elements can be expressed as a linear combination of those two basis elements, with coefficients in a much smaller subfield. \\
\end{abstract}
}

\begin{small}
\textcolor{secheader}{
\noindent {\em{Keywords:}}~~ Finite fields; ~Algebraic field extensions; ~Normal basis; ~Optimal normal basis. 
}
\end{small}

\textcolor{titlecolor}{\section{\bf{Introduction}}}
In this paper, a method for constructing a near optimal normal basis for an algebraic extension of specific dimension (degree) over a finite field is described. The optimality criteria are that the multiplication tables have as few nonzero entries as possible. The extensions can be classified as either Artin-Schreier extensions, where the degree of extension is the same the characteristic of the field, or other extensions, where the degree of extension is relatively prime with the characteristic of the field. An application to discrete logarithm cryptography is also described, in a later section. Previous results on optimal normal bases for finite field extensions are mostly based on those studied in \cite{MOVW1988}.  Algorithms for construction of finite fields of specified number elements are described in \cite{AL1986, Shoup1990}, and randomized algorithms in \cite{Rabin1980, Shoup1994}, while permutation polynomials and irreducible polynomials over a given finite field are presented in \cite{GM1994}.

\textcolor{titlecolor}{\section{\label{Sec-Artin-Schreier-Extension}\bf{Artin-Schreier Extensions of Finite Fields}}}
Throughout the paper, let $\fieldChar$ be a fixed prime number as well as the characteristic of a finite field
$\scalars$, and $\basel{\Integers}{\fieldChar}$ be the prime field with elements $0,\, \ldots, \, \fieldChar-2$
and $\fieldChar-1$, equipped with the arithmetic operations of addition and multiplication modulo $\fieldChar$.

\textcolor{titlecolor}{
\subsection{\label{Sec-Normality}Normality of the Bases Generated by Artin-Schreier Extensions} 
}
The lemma below plays an important role in the results that follow:

\textcolor{secheader}{
\begin{lemma}\label{Lemma-1}
Let $\scalars$ and $\extensionField$ be finite fields containing
$\fieldChar^{n}$ and $\fieldChar^{mn}$ elements, respectively,
for some prime number $\fieldChar$ and positive integers $m$ and $n$,
where $m \geq 2$. Let $\{\delta^{\fieldChar^{i\cdot n}}\,:\, 0 \leq i \leq m-1\}$ 
be a basis for the extension field $\extensionField$ over $\scalars$,
for some $\delta \in \extensionField$. Then, for every $d \in \scalars$,
such that $\bgls md + \bglb\sum_{j = 0}^{m-1} \delta^{\fieldChar^{j\cdot n}}\bgrb \bgrs \neq 0$,
the set $\{(\delta+d)^{\fieldChar^{i\cdot n}}\,:\, 0 \leq i \leq m-1\}$
is also a basis for $\extensionField$ as an extension field of $\scalars$.
\end{lemma}
}

\noindent \proof  The linear span of the set 
$\{(\delta+d)^{\fieldChar^{i\cdot n}}\, : \, 0 \leq i \leq m-1\}$
is the same as that of 
$\{\delta+d\} \cup \{(\delta+d)^{\fieldChar^{i\cdot n}}-(\delta+d)^{\fieldChar^{(i-1)\cdot n}}\, : \, 1 \leq i \leq m-1\}$,
which, in turn, is that of 
$\{\delta+d\} \cup \{\delta^{\fieldChar^{i\cdot n}}-\delta^{\fieldChar^{(i-1)\cdot n}}\, : \, 1 \leq i \leq m-1\}$. 
When $d = 0$, the set $\{\delta^{\fieldChar^{i\cdot n}}-\delta^{\fieldChar^{(i-1)\cdot n}}\, : \, 1 \leq i \leq m-1\}$ is linearly independent over $\scalars$, by the hypothesis.
Now, for some sequence of elements $ \basel{c}{i} \in \scalars$, $1 \leq i \leq m-1$,
if $(\delta+d)$
$  = $
$ \sum_{i = 1}^{m-1} \basel{c}{i}
 \bglb  \delta^{\fieldChar^{i\cdot n}}- \delta^{\fieldChar^{(i-1)\cdot n}}\bgrb $,
 then
\begin{small}
\begin{eqnarray*}
&& d ~~ = ~~ -(1+\basel{c}{1})\delta ~ + ~ \sum_{i = 1}^{m-2}(\basel{c}{i} - \basel{c}{i+1}) \delta^{\fieldChar^{i\cdot n}} ~ + ~ \basel{c}{m-1}\delta^{\fieldChar^{(m-1)\cdot n}}
 \end{eqnarray*}
\end{small} 
\lspace However, since
$d = d \bglb\sum_{j = 0}^{m-1} \delta^{\fieldChar^{j\cdot n}}\bgrb^{-1} \times $
$\sum_{i = 0}^{m-1} \delta^{\fieldChar^{i\cdot n}}$ is the unique expression for $d$,
as a linear combination of $\delta^{\fieldChar^{i\cdot n}}$, for $0 \leq i \leq m-1$,
it follows that 
\begin{small}
\begin{eqnarray*}
&&\ltab  \basel{c}{1} ~ = ~ -(\tau+1)\, , ~~ \basel{c}{i+1} ~ = ~ \basel{c}{i}- \tau\,, ~~ \textrm{for}~~ 1 \leq i \leq m-2\,, ~~  \textrm{and}~~ \basel{c}{m-1} ~ = ~ \tau,
 \end{eqnarray*}
\end{small} 
\lspace where $\tau = d \bglb\sum_{j = 0}^{m-1} \delta^{\fieldChar^{j\cdot n}}\bgrb^{-1}$. 
By induction on $i$, it can be deduced from the first two requirements that
$\basel{c}{i} = -(1+i\cdot \tau)$, for $1 \leq i \leq m-1$.
Now, for the last requirement that $\basel{c}{m-1} = \tau$
to be consistent with that $\basel{c}{m-1}= -(1+(m-1)\tau)$,
as deduced from the previous requirements,
it needs to be true that $m\tau = -1$,
contrary to the hypothesis on $m$ and $d$.
Thus, the coefficients $\basel{c}{i} \in \scalars$,
for $1 \leq i \leq m-1$, such that $(\delta+d)$
$  = $
$ \sum_{i = 1}^{\fieldChar-1} \basel{c}{i}
 \bglb  \delta^{\fieldChar^{i\cdot n}}- \delta^{\fieldChar^{(i-1)\cdot n}}\bgrb $,
cannot exist, assuming that the set
$\{\delta^{\fieldChar^{i\cdot n}}\, : \, 0 \leq i \leq m-1\}$
is linearly independent over $\scalars$. \qed
\\

The following result allows construction of many normal elements in an Artin-Schreier extension of a finite field:

\textcolor{secheader}{
\begin{theorem}
\label{Normal-Basis-of-Sec-Artin-Schreier-Extension-of-Finite-Fields-1}
{~~\bf{{\sf(Normal Basis Generated by the Artin-Schreier-Extension)}}~~}
Let $\scalars$ be a finite dimensional extension field
of $\basel{\Integers}{\fieldChar}$, of vector space
dimension $n$, for some positive integer $n$. Let
$\alpha,\, \varepsilon \in \scalars\backslash \{0\}$ be
such that the polynomial $x^{\fieldChar}-x-\alpha$ is
irreducible in $\singlevariablepolynomials{\scalars}{x}$,
and the set  $\left\{\varepsilon^{\fieldChar^{j}}\,:\, 0 \leq j \leq n-1\right\}$
is a basis for $\scalars$ as an extension field of $\basel{\Integers}{\fieldChar}$.
Let  $\extensionField = \singlevariablepolynomials{\scalars}{\beta}
/\bglb \beta^{\fieldChar}-\beta-\alpha\bgrb$.
The following statements hold:
\begin{enumerate}
\item for every $c \in \scalars$, the set
$\left\{\bglb\beta^{\fieldChar-1}-c\bgrb^{\fieldChar^{i\cdot n}}\,:\,
0 \leq i \leq \fieldChar-1\right\}$ is a basis for $\extensionField$
as an extension field of $\scalars$, ~~and
\item  for every $c \in \basel{\Integers}{\fieldChar}$,  the set
$\left\{\bglb\varepsilon(\beta^{\fieldChar-1}-c)\bgrb^{\fieldChar^{l}}\,:\,
0 \leq l \leq n\fieldChar-1\right\}$ is a basis for $\extensionField$
as an extension field of $\basel{\Integers}{\fieldChar}$.
\end{enumerate}
\end{theorem}
}

\noindent \proof {\underline{\em{Part (1)}}~~} The linear span of  the set 
$\left\{\bglb\beta^{\fieldChar-1}\bgrb^{\fieldChar^{i\cdot n}}\,:\,
0 \leq i \leq \fieldChar-1\right\}$ with coefficients in $\scalars$
is the same as that of the set
$\{\beta^{\fieldChar-1}\} ~\bigcup~
\{\bglb\beta^{\fieldChar-1}\bgrb^{\fieldChar^{i\cdot n}} - 
\bglb\beta^{\fieldChar-1}\bgrb\,:\,
1 \leq i \leq \fieldChar-1\}$.
 Now, $\beta^{\fieldChar^{l}}  =
\beta +  \sum_{i = 0}^{l-1}\alpha^{\fieldChar^{i}}$,
for $0 \leq l \leq n\fieldChar-1$. Let
$\basel{s}{j} = \sum_{i = 0}^{j-1} \alpha^{\fieldChar^{i}}$,
$0 \leq j \leq n-1$, and $h =  \sum_{i = 0}^{n-1} \alpha^{\fieldChar^{i}}$.
Since $x^{\fieldChar}-x-\alpha$ is irreducible in
$\singlevariablepolynomials{\scalars}{x}$, it follows that
$h \in \basel{\Integers}{\fieldChar} \backslash \{0\}$ 
(Theorem 3.78 and Corollary 3.79 in \cite{LN1986}).
Thus, $\basel{h}{i\cdot n + j} = \basel{s}{j}+i\cdot h$, and
$\beta^{\fieldChar^{j+i\cdot n}}  = \beta+ \basel{s}{j} + i \cdot h$, 
for $0 \leq i \leq \fieldChar-1$ and $0 \leq j \leq n-1$.
Now,  $\beta^{(\fieldChar-1)\fieldChar^{i\cdot n}}  = 
\bglb \beta + i\cdot h\bgrb^{\fieldChar-1}  =  
\beta^{\fieldChar-1} + \sum_{k = 1}^{\fieldChar-1} \frac{(\fieldChar-1)!}{k!(\fieldChar-k-1)!} 
(i\cdot h)^{k} \beta^{\fieldChar-k-1}$, for $1 \leq i \leq \fieldChar-1$.
Thus, $\bglb\beta^{(\fieldChar-1)\fieldChar^{i\cdot n}} - \beta^{\fieldChar-1}\bgrb
 =   \sum_{k = 1}^{\fieldChar-1} \frac{(\fieldChar-1)!}{k!(\fieldChar-k-1)!} 
(i\cdot h)^{k} \beta^{\fieldChar-k-1}$, for $1 \leq i \leq \fieldChar-1$.
Now,  for $1 \leq i,\, k \leq \fieldChar-1$,
$\frac{(\fieldChar-1)!}{k!(\fieldChar-k-1)!} \neq 0$,  and
the matrix with $(i\cdot h)^{k}$ as the element in the $i$-th row
and $k$-th column can be expressed as the product of a permutation
matrix and a Vandermonde matrix, and hence the corresponding
linear transformation becomes invertible. Thus, the
linear span with coefficients in $\scalars$ (and respectively
also with coefficients in $\basel{\Integers}{\fieldChar}$) of
the set $\{\bglb\beta^{\fieldChar-1}\bgrb^{\fieldChar^{i\cdot n}} - 
\bglb\beta^{\fieldChar-1}\bgrb\,:\, 1 \leq i \leq \fieldChar-1\}$
is the same as that of the set
$\{\beta^{i-1}\,:\,\, 1 \leq i \leq \fieldChar-1\}$, and the set 
$\{\bglb\beta^{(\fieldChar-1)}\bgrb^{\fieldChar^{i\cdot n}}\,:\,
0 \leq i \leq \fieldChar-1\}$ is a basis for $\extensionField$
 as an extension field of $\scalars$, since
the set $\{\beta^{i-1}\,:\,\, 1 \leq i \leq \fieldChar-1\}
 \cup \{\beta^{\fieldChar-1}\}$ is so. 
 The result follows by applying Lemma \ref{Lemma-1} with 
$\delta = \beta^{\fieldChar-1}$, $~d = -c~$ and $~m = \fieldChar$. 
\\
\noindent {\underline{\em{Part (2)}}~~}
By the previous part, 
and for any $c \in \basel{\Integers}{\fieldChar}$,
 the linear span with coefficients in
 $\basel{\Integers}{\fieldChar}$ of the set
$\{\bglb \beta^{\fieldChar-1} - c  \bgrb^{\fieldChar^{i\cdot n}}\,:\,
1 \leq i \leq \fieldChar-1\}$ is the same as that of the set
$\{\bglb\beta^{i-1}\bgrb\,:\, 1 \leq i \leq \fieldChar-1\}$.
Now, since $1 = \beta^{0}$ belongs to the latter, it follows that,
for any $c \in \basel{\Integers}{\fieldChar}$,
 the linear span
with coefficients in $\basel{\Integers}{\fieldChar}$ of the set
$\{\bglb \beta^{\fieldChar-1} - c  \bgrb^{\fieldChar^{i\cdot n}}\,:\,
0 \leq i \leq \fieldChar-1\}$ is the same as that of the set
$\{\bglb\beta^{i}\bgrb\,:\, 0 \leq i \leq \fieldChar-1\}$.
Therefore, (A) for every $\varepsilon \in \extensionField$,
$c \in \basel{\Integers}{\fieldChar}$ and integer $j$, such that
$0 \leq j \leq n-1$, the linear span with coefficients in
 $\basel{\Integers}{\fieldChar}$ of the set
$\{\bglb \varepsilon(\beta^{\fieldChar-1}-c)^{\fieldChar^{i\cdot n}}\bgrb^{\fieldChar^{j}}\,:\,
0 \leq i \leq \fieldChar-1\}$ is the same as that of the set
$\{\bglb\varepsilon\beta^{i}\bgrb^{\fieldChar^{j}}\,:\, 0 \leq i \leq \fieldChar-1\}$,
and hence, (B) the linear span with coefficients in $\basel{\Integers}{\fieldChar}$
of the set 
 $\{\bglb \varepsilon(\beta^{\fieldChar-1}-c)^{\fieldChar^{i\cdot n}}\bgrb^{\fieldChar^{j}}\,:\,
0 \leq i \leq \fieldChar-1,\, 0 \leq j \leq n-1\}$ is the same as
that of the set $\{\bglb\varepsilon\beta^{i}\bgrb^{\fieldChar^{j}}
\,:\, 0 \leq i \leq \fieldChar-1,\, 0 \leq j \leq n-1\}$.
 For $\basel{a}{i,\, j} \in \basel{\Integers}{\fieldChar}$,
  $0 \leq i \leq \fieldChar-1$,  $0 \leq j \leq n-1$, 
  \begin{small}
  \begin{eqnarray*}
&&\ltab \ltab  \sum_{i = 0}^{\fieldChar-1} \sum_{j = 0}^{n-1} \basel{a}{i,\, j}
\bglb \varepsilon \beta^{i} \bgrb^{\fieldChar^{j}}  ~~ = ~~
  \sum_{j = 0}^{n-1} \basel{a}{0,\, j} \varepsilon^{\fieldChar^{j}} ~ + ~
 \sum_{i = 1}^{\fieldChar-1} \sum_{j = 0}^{n-1} \basel{a}{i,\, j}
\bglb \varepsilon \beta^{i} \bgrb^{\fieldChar^{j}}\\
&& = ~~
\sum_{j = 0}^{n-1} \basel{a}{0,\, j} \varepsilon^{\fieldChar^{j}} ~ + ~
 \sum_{i = 1}^{\fieldChar-1} \sum_{j = 0}^{n-1} \basel{a}{i,\, j} 
 \varepsilon^{\fieldChar^{j}}\bglb \beta+\basel{s}{j}\bgrb^{i}\,, \tab \textrm{where}~~ \basel{s}{j} = \sum_{l = 0}^{j-1}\alpha^{\fieldChar^{l}} \\
 &&  = ~~
\sum_{j = 0}^{n-1} \basel{a}{0,\, j} \varepsilon^{\fieldChar^{j}} ~ + ~
 \sum_{i = 1}^{\fieldChar-1} \sum_{j = 0}^{n-1} \basel{a}{i,\, j}
  \varepsilon^{\fieldChar^{j}} 
 \sum_{k = 0}^{i}\frac{i!}{k!(i-k)!} \basel{s^{i-k}}{j} \beta^{k}\\
&& = ~~
\sum_{j = 0}^{n-1} \basel{a}{0,\, j} \varepsilon^{\fieldChar^{j}} ~ + ~
 \sum_{i = 1}^{\fieldChar-1} \sum_{j = 0}^{n-1}
 \basel{a}{i,\, j} \varepsilon^{\fieldChar^{j}} \beta^{i} ~ + ~
\sum_{i = 1}^{\fieldChar-1} \sum_{j = 0}^{n-1}
  \sum_{k = 0}^{i-1} \frac{i!~ \basel{a}{i,\, j}}{k!(i-k)!}
   \basel{s^{i-k}}{j} \varepsilon^{\fieldChar^{j}} \beta^{k}
   \end{eqnarray*}
  \end{small}
\lspace   Now,  $ \sum_{i = 1}^{\fieldChar-1} \sum_{j = 0}^{n-1}
 \basel{a}{i,\, j} \varepsilon^{\fieldChar^{j}} \beta^{i} ~~ = ~~
 \sum_{j = 0}^{n-1} \basel{a}{\fieldChar-1,\, j}
  \varepsilon^{\fieldChar^{j}} \beta^{\fieldChar-1} ~+~
 \sum_{i = 1}^{\fieldChar-2} \sum_{j = 0}^{n-1}
 \basel{a}{i,\, j} \varepsilon^{\fieldChar^{j}} \beta^{i}$,
 and 
\begin{small}
\begin{eqnarray*}
&& \shiftleft \ltab \ltab  \sum_{i = 1}^{\fieldChar-1} \sum_{j = 0}^{n-1}
  \sum_{k = 0}^{i-1} \frac{i!~ \basel{a}{i,\, j}}{k!(i-k)!}
   \basel{s^{i-k}}{j} \varepsilon^{\fieldChar^{j}} \beta^{k}
     ~~ = ~~ 
     \sum_{k = 0}^{\fieldChar-2}  \beta^{k} \sum_{j = 0}^{n-1} 
   \varepsilon^{\fieldChar^{j}}  \sum_{i = k+1}^{\fieldChar-1}
       \frac{i!~ \basel{a}{i,\, j}}{k!(i-k)!} \basel{s^{i-k}}{j}\\
&&  = ~~~~ \sum_{j = 0}^{n-1} 
   \varepsilon^{\fieldChar^{j}}  \sum_{i = 1}^{\fieldChar-1}
        \basel{a}{i,\, j} \basel{s^{i}}{j} ~~ + ~~
        \sum_{k = 1}^{\fieldChar-2}  \beta^{k} \sum_{j = 0}^{n-1} 
   \varepsilon^{\fieldChar^{j}}  \sum_{i = k+1}^{\fieldChar-1}
       \frac{i!~ \basel{a}{i,\, j}}{k!(i-k)!} \basel{s^{i-k}}{j}
\end{eqnarray*}
\end{small}
\lspace  Thus
  \begin{small}
 \begin{eqnarray*}
  && \ltab\ltab~~
  \sum_{i = 0}^{\fieldChar-1} \sum_{j = 0}^{n-1} \basel{a}{i,\, j}
\bglb \varepsilon \beta^{i} \bgrb^{\fieldChar^{j}}  ~~ = ~~
\sum_{j = 0}^{n-1} \basel{a}{0,\, j} \varepsilon^{\fieldChar^{j}} ~+~
 \sum_{j = 0}^{n-1} \basel{a}{\fieldChar-1,\, j} 
 \varepsilon^{\fieldChar^{j}} \beta^{\fieldChar-1}  ~+~
 \sum_{k = 1}^{\fieldChar-2} \sum_{j = 0}^{n-1} \basel{a}{k,\, j}
 \varepsilon^{\fieldChar^{j}} \beta^{k} \\
 && \shiftright ~~~~ +~ \sum_{j = 0}^{n-1} \sum_{i = 1}^{\fieldChar-1}
 \basel{a}{i,\, j}  \varepsilon^{\fieldChar^{j}} \basel{s^{i}}{j} 
  ~ + ~ \sum_{k = 1}^{\fieldChar-2}   \beta^{k} \sum_{j = 0}^{n-1}
 \varepsilon^{\fieldChar^{j}} \sum_{i = k+1}^{\fieldChar-1}
 \frac{i!~ \basel{a}{i,\, j}}{k!(i-k)!} \basel{s^{i-k}}{j}\\
 && \tab = ~~ \sum_{j = 0}^{n-1} \basel{a}{\fieldChar-1,\, j} 
 \varepsilon^{\fieldChar^{j}} \beta^{\fieldChar-1}  ~ + ~
 \sum_{k = 1}^{\fieldChar-2} \beta^{k} \sum_{j = 0}^{n-1}
  \varepsilon^{\fieldChar^{j}} \bglb \basel{a}{k,\, j}  
+  \sum_{i = k+1}^{\fieldChar-1}
  \frac{i!~ \basel{a}{i,\, j}}{k!(i-k)!} \basel{s^{i-k}}{j}\bgrb\\
&& \shiftright \shiftright + ~ \sum_{j = 0}^{n-1} \sum_{i = 1}^{\fieldChar-1}
 \basel{c}{i,\, j}  \varepsilon^{\fieldChar^{j}} \basel{s^{i}}{j} +
\sum_{j = 0}^{n-1} \basel{a}{0,\, j} \varepsilon^{\fieldChar^{j}}
\end{eqnarray*}
\end{small}
\lspace Now, the set $ \{\varepsilon^{\fieldChar^{j}}\beta^{k}\,:\,
0 \leq j \leq n-1,\, 0\leq k \leq \fieldChar-1 \}$
is a basis for $\extensionField$ with $\basel{\Integers}{\fieldChar}$
as the field. If the above expression evaluates to $0$, then,
comparing the coefficients of $\varepsilon^{\fieldChar^{j}}\beta^{\fieldChar-1-l}$,
for $l = 0,\, 1,\,\ldots,\, \fieldChar-2$, it follows that
$\basel{a}{\fieldChar-1,\, j} = 0$, for $0 \leq j \leq n-1$,
and subsequently, if $\basel{a}{i,\,j} = 0$,
for $k+1 \leq i \leq \fieldChar-1$, for some index $k$,
occurring in the lower bound of the running index $i$,
where $1 \leq k \leq \fieldChar-2$, then $\basel{a}{k,\,j} = 0$,
for $0 \leq j \leq n-1$. Now, since $\basel{a}{k,\, j} = 0$,
for $1 \leq k \leq \fieldChar-1$ and $0 \leq j \leq n-1$,
$\sum_{j = 0}^{n-1} \sum_{i = 1}^{\fieldChar-1}
\basel{a}{i,\, j} \varepsilon^{\fieldChar^{j}}\basel{s^{i}}{j} = 0$,
and hence $\basel{a}{0,\,j} = 0$, for $0 \leq j \leq n-1$.
Thus, the set $\{\bglb\varepsilon\beta^{i}\bgrb^{\fieldChar^{j}}
\,:\, 0 \leq i \leq \fieldChar-1,\, 0 \leq j \leq n-1 \}$,
and hence the set 
$\{\bglb \varepsilon(\beta^{\fieldChar-1}-c)^{\fieldChar^{i\cdot n}}\bgrb^{\fieldChar^{j}}\,:\,
0 \leq i \leq \fieldChar-1,\, 0 \leq j \leq n-1\}$,
is a basis for $\extensionField$ as an extension field of
$\basel{\Integers}{\fieldChar}$.\qed
\\

The following result is crucial for inductive construction of
multiple  Artin-Schreier extensions, preserving both normality
and optimality.

\textcolor{secheader} {
\begin{theorem}
\label{Normal-Basis-of-Sec-Artin-Schreier-Extension-of-Finite-Fields-2}
  {\bf {\textsf{~~(Inductive Construction of Normal Bases of Multiple Artin-Schreier-Extensions)~~}}} 
Let $\scalars$ be a finite dimensional extension field
of $\basel{\Integers}{\fieldChar}$, of vector space
dimension $n$, such that $n = m\fieldChar$, for some
positive integer $m$. Let $\alpha \in \scalars\backslash \{0\}$
be such that the polynomial $x^{\fieldChar}-x-\alpha$
is irreducible in $\singlevariablepolynomials{\scalars}{x}$,
and the set $\{\alpha^{-\fieldChar^{i}}\,:\, 0 \leq i \leq n-1\}$
is a basis for $\scalars$ as an extension field of
$\basel{\Integers}{\fieldChar}$. Let
$\extensionField = \singlevariablepolynomials{\scalars}{\beta}
/\bglb \beta^{\fieldChar}-\beta-\alpha\bgrb$, and
$\delta = \bglb\beta^{-1}-b\bgrb^{-1}$, for some
$b \in \basel{\Integers}{\fieldChar} \backslash \{0\}$.
Then, the polynomial $x^{\fieldChar}-x-\delta$ is irreducible
in $\singlevariablepolynomials{\extensionField}{x}$, and 
the set $\{\delta^{-\fieldChar^{i}}\,:\, 0 \leq i \leq n\fieldChar-1\}$
is a basis for $\extensionField$ as an extension field of
$\basel{\Integers}{\fieldChar}$.
\end{theorem}
}
\proof By Theorem \ref{Normal-Basis-of-Sec-Artin-Schreier-Extension-of-Finite-Fields-1},
the element $\alpha^{-1}\bglb\beta^{\fieldChar-1}-1\bgrb = 
\alpha^{-1}\bglb\alpha\beta^{-1}\bgrb = \beta^{-1}$ generates
a normal basis for $\extensionField$ as an extension field of
$\basel{\Integers}{\fieldChar}$. Now, by Lemma \ref{Lemma-1},
by appropriate reinterpretation of parameters in the statements,
$\delta^{-1}$ generates a normal basis for $\extensionField$
as an extension field of $\basel{\Integers}{\fieldChar}$. 
The polynomial $x^{\fieldChar}-x-\delta$ is irreducible in
$\singlevariablepolynomials{\extensionField}{x}$ if and only if
the coefficient of $x^{n\fieldChar-1}$ in the minimal polynomial
of $\delta$ in $\singlevariablepolynomials{\basel{\Integers}{\fieldChar}}{x}$
does not vanish \cite{LN1986}.
Let $\basel{m}{t}(x)$ be the minimal polynomial of $t$ in
$\singlevariablepolynomials{\basel{\Integers}{\fieldChar}}{x}$,
for $t \in \extensionField$. Now, $\basel{m}{\delta}(x) = 
\bglb\basel{m}{\delta^{-1}}(0)\bgrb^{-1} x^{n\fieldChar}\basel{m}{\delta^{-1}}\bglb x^{-1}\bgrb$,
and the coefficient of $x^{n\fieldChar-1}$ in $\basel{m}{\delta}(x)$
is nonzero  if and only if  the coefficient of $x$ in
$\basel{m}{\delta^{-1}}(x)$ is nonzero.
Let $\basel{m}{\alpha}(x) = x^{n}+\sum_{i = 1}^{n}\basel{a}{n-i} x^{n-i}$,
for some $\basel{a}{i} \in \basel{\Integers}{\fieldChar}$, $0 \leq i \leq n-1$.
Now, $\basel{a}{0} \neq 0$, $\basel{a}{1} \neq 0$
and $\basel{a}{n-1} \neq 0$ in $\basel{m}{\alpha}(x)$\,,~
 $\basel{m}{\delta^{-1}}(x) = \basel{m}{\beta^{-1}}(x+b)$\,,~
 $\basel{m}{\beta^{-1}}(x) =  \bglb \basel{m}{\beta}(0)\bgrb^{-1}x^{n\fieldChar}\basel{m}{\beta}(x^{-1})$\,,~
$\basel{m}{\beta}(x)  =  \basel{m}{\alpha}\bglb x^{\fieldChar}-x \bgrb
 = \bglb x^{\fieldChar}-x\bgrb^{n}+
\sum_{i = 1}^{n}\basel{a}{n-i} \bglb x^{\fieldChar}-x\bgrb^{n-i}$\,,~
 $\basel{m}{\beta}(0)  =  \basel{a}{0}$\,,~ and 
\begin{small}
\begin{eqnarray*}
&& \ltab\ltab ~~~~
\basel{a}{0}\basel{m}{\beta^{-1}}(x) ~~ = ~~
\basel{m}{\beta}(0) \basel{m}{\beta^{-1}}(x) ~~ = ~~
x^{n\fieldChar}\basel{m}{\beta}(x^{-1}) \\
&&  = ~~
x^{n\fieldChar} \bgls
\bglb x^{-\fieldChar}-x^{-1}\bgrb^{n}+\sum_{i = 1}^{n}\basel{a}{n-i} \bglb x^{-\fieldChar}-x^{-1}\bgrb^{n-i}\bgrs \\
&&  = ~~ \bglb 1-x^{\fieldChar-1}\bgrb^{n}+\sum_{i = 1}^{n}\basel{a}{n-i}
x^{i\fieldChar} \bglb 1-x^{\fieldChar-1}\bgrb^{n-i} ~~~~ = ~~~~
\bglb 1-x^{\fieldChar-1}\bgrb^{n} ~~+ \\
&& \ltab \sum_{i = 1}^{n-2}\basel{a}{n-i}
x^{i\fieldChar-2} \bglb x - x^{\fieldChar}\bgrb^{2} \bglb 1-x^{\fieldChar-1}\bgrb^{n-i-2}
+ \basel{a}{1} x^{(n-1)\fieldChar-1} \bglb x - x^{\fieldChar}\bgrb ~ + ~ \basel{a}{0} x^{n\fieldChar}\,,
 ~~ \textrm{and}\\
&& \ltab\ltab ~~~~
\basel{a}{0}\basel{m}{\delta^{-1}}(x) ~~ = ~~
\basel{m}{\beta}(0) \basel{m}{\beta^{-1}}(x+b)\\
&& = ~~
\sum_{i = 1}^{n-2}\basel{a}{n-i}
(x+b)^{i\fieldChar-2} \bglb x - x^{\fieldChar}\bgrb^{2} \bglb 1-(x+b)^{\fieldChar-1}\bgrb^{n-i-2} 
~~~~+\\
&& \tab \tab \bglb 1-(x+b)^{\fieldChar-1}\bgrb^{n} 
 ~~+ ~~\basel{a}{1} (x+b)^{(n-1)\fieldChar-1} \bglb x - x^{\fieldChar}\bgrb
~~+~~ \basel{a}{0} (x+b)^{n\fieldChar}\,, 
\end{eqnarray*}
\end{small}
\lspace 
since $(x+b)^{\fieldChar}-(x+b) = x^{\fieldChar}-x$. Now, 
$(x+b)^{i\fieldChar-2} \bglb x - x^{\fieldChar}\bgrb^{2}
\bglb 1-(x+b)^{\fieldChar-1}\bgrb^{n-i-2}$ is divisible
by $x^{2}$, for $1 \leq i \leq n-2$,~
$\bglb 1-(x+b)^{\fieldChar-1}\bgrb^{n}  =  
\bglb 1-(x^{\fieldChar}+b)^{\fieldChar-1}\bgrb^{m}~$,
since $n = m\fieldChar$,~ and ~
 $(x+b)^{n\fieldChar} =  (x^{\fieldChar}+b)^{n}~$.
Thus, the coefficient of $x$ in each of the 
polynomials $\bglb 1-(x+b)^{\fieldChar-1}\bgrb^{n}$\,,~ 
$(x+b)^{n\fieldChar}~$ ~ and~ 
$(x+b)^{i\fieldChar-2} \bglb x - x^{\fieldChar}\bgrb^{2}
\bglb 1-(x+b)^{\fieldChar-1}\bgrb^{n-i-2}$,
for $1 \leq i \leq n-2$, is $0$. However,
the coefficient of $x$ in
$\basel{a}{1} (x+b)^{(n-1)\fieldChar-1} \bglb x - x^{\fieldChar}\bgrb $
is $\basel{a}{1}b^{(n-1)\fieldChar-1} \neq 0$,
since $\basel{a}{1} \neq 0$ and
$b \in \basel{\Integers}{\fieldChar} \backslash \{0\}$.
Now as $\sum_{i = 0}^{n\fieldChar-1} \delta^{\fieldChar^{i}} \neq 0$,
the polynomial $x^{\fieldChar}-x-\delta \in \singlevariablepolynomials{\extensionField}{x}$
is irreducible, such that $\delta^{-1}$ generates a normal basis for $\extensionField$
as an extension field of $\basel{\Integers}{\fieldChar}$. \qed

\textcolor{titlecolor}{
\subsection{\label{Sec-optimality}Near Optimality of the Normal Bases of Artin-Schreier Extensions}
}

By Theorem \ref{Normal-Basis-of-Sec-Artin-Schreier-Extension-of-Finite-Fields-1} discussed
in the preceding subsection, with $c = 1$, the element $\alpha \beta^{-1}$ is a normal element in $\extensionField = \singlevariablepolynomials{\scalars}{\beta}/(\beta^{\fieldChar}-\beta-\alpha)$ over $\scalars$, in the notation followed there.  Thus, the element $\gamma = \beta^{-1}$,  satisfying the equation $-\alpha^{-1}\gamma^{\fieldChar}\bglb \gamma^{-\fieldChar}-\gamma^{-1}-\alpha\bgrb = 0$, is a normal element in $\extensionField$.
The element $\beta$ also satisfies the equation $\beta^{\fieldChar^{n}}-\beta-h = 0$, where $h = \sum_{i = 0}^{n-1} \alpha^{\fieldChar^{i}}$, and therefore, for $1 \leq j \leq \fieldChar-1$, the element $\beta$ satisfies the equation $\beta^{\fieldChar^{j\cdot n}}-\beta^{\fieldChar^{(j-1)\cdot n}} - h = 0$, and summing for $1 \leq j \leq i$, for some integer $i$, such that $1 \leq i \leq \fieldChar-1$, it follows that $\beta^{\fieldChar^{i\cdot n}}-\beta - i\cdot h = 0$, and therefore, $\gamma$ satisfies the equation $\gamma^{-\fieldChar^{i\cdot n}}-\gamma^{-1} - i\cdot h = 0$, and multiplying throughout the last equation by  $(i\cdot h)^{-1}\gamma^{\fieldChar^{i\cdot n}+1}$,  the following is obtained:
\begin{small}
\begin{equation}
\gamma^{\fieldChar^{i\cdot n}+1}~~ = ~~ (i\cdot h)^{-1}\gamma-(i\cdot h)^{-1}\gamma^{\fieldChar^{i\cdot n}} 
\,, \tab \textrm{for}~~ 1 \leq i \leq \fieldChar-1 \label{Non-square-product-of-normal-basis-elements}
\end{equation}
\end{small}
\lspace As for $\gamma^{2}$, from the equations that $\bglb \sum_{i = 0}^{\fieldChar-1}\gamma^{\fieldChar^{i\cdot n}} \bgrb \gamma$
$ = $
$\sum_{i = 0}^{\fieldChar-1}\gamma^{\fieldChar^{i\cdot n}+1}$ 
$ = $
$\gamma^{2} + 
\sum_{i = 1}^{\fieldChar-1}(i\cdot h)^{-1}\bgls \gamma-\gamma^{\fieldChar^{i\cdot n}} \bgrs$,  
it follows that 
\begin{small}
\begin{eqnarray*}
&& \ltab  \gamma^{2} ~~ = ~~
  \bgls \sum_{i = 0}^{\fieldChar-1}\gamma^{\fieldChar^{i\cdot n}}   ~ - ~
 \sum_{i = 1}^{\fieldChar-1} (i\cdot h)^{-1} \bgrs \gamma ~ + ~  \sum_{i = 1}^{\fieldChar-1} (i\cdot h)^{-1}\gamma^{\fieldChar^{i\cdot n}} \\
&&  = ~~  \bgls \sum_{i = 0}^{\fieldChar-1}\gamma^{\fieldChar^{i\cdot n}}   ~ - ~
 h^{-1} \cdot \sum_{i = 1}^{\fieldChar-1} i \bgrs \gamma ~ + ~  \sum_{i = 1}^{\fieldChar-1} (i\cdot h)^{-1}\gamma^{\fieldChar^{i\cdot n}} \\
&&   = ~~  \left  \{ \begin{array}{l}  
 \bgls \sum_{i = 0}^{\fieldChar-1}\gamma^{\fieldChar^{i\cdot n}}   ~ - ~
 h^{-1} \bgrs \gamma ~ + ~  h^{-1}\cdot \gamma^{\fieldChar^{n}}\,, \tab \textrm{when}~~ \fieldChar = 2\,, \tab \textrm{and} \\
 \bgls \sum_{i = 0}^{\fieldChar-1}\gamma^{\fieldChar^{i\cdot n}} \bgrs \gamma ~ + ~  \sum_{i = 1}^{\fieldChar-1} (i\cdot h)^{-1}\gamma^{\fieldChar^{i\cdot n}}\,, \tab \textrm{when}~~ \fieldChar \geq 3 
 \end{array} \right .
 \end{eqnarray*}
 \end{small}
 \lspace Since $\gamma^{\fieldChar} +\alpha^{-1} \gamma^{\fieldChar-1} -\alpha^{-1} = 0$, it follows that $\sum_{i = 0}^{\fieldChar-1} \gamma^{\fieldChar^{i\cdot n}} = -\alpha^{-1} = -\prod_{i = 0}^{\fieldChar-1} \gamma^{\fieldChar^{i\cdot n}}$, and the remaining coefficients in the product are all elements in the prime field $\basel{\Integers}{\fieldChar}$. Now, if $n = \fieldChar t$, for some positive integer $t$, and $\fieldChar \geq 3$, then the normal element $\alpha^{-1}$ can be chosen such that the normal basis generated by it for $\scalars$, over its subfield of $\fieldChar^{t}$ elements, is nearly optimal, similar to $\gamma$ for $\extensionField$. However, some more little work needs to be done, with minor adjustments for preserving optimality of the  multiplication tables as nearly as possible, for inductive application.
 
Let $\delta = (\gamma-b)^{-1}$, for some $b \in \scalars$.
The element $\delta^{-1}$ is a normal element in $\extensionField$
as an extension field of $\scalars$, for any $b \in \scalars$,
by Lemma \ref{Lemma-1} of the preceding subsection. It is
convenient and preferable to choose $b = 1$, and, for applying
Theorem \ref{Normal-Basis-of-Sec-Artin-Schreier-Extension-of-Finite-Fields-2},
 $b$ must be chosen to be an integer in  $\finiteIntegerField \backslash \{0\}$.
 Now,
for $1 \leq i \leq \fieldChar-1$,
\begin{eqnarray}
&& \ltab \ltab \ltab \ltab  \delta^{-(1+\fieldChar^{i\cdot n})} ~~ = ~~
(\gamma-b)(\gamma^{\fieldChar^{i\cdot n}}-b) ~~ = ~~
 \gamma^{(1+\fieldChar^{i\cdot n})} - b \bglb \gamma + \gamma^{\fieldChar^{i\cdot n}}\bgrb + b^{2}  \nonumber   \\
&&  = ~~ \bglb (i\cdot h)^{-1} - b \bgrb \gamma - \bglb (i\cdot h)^{-1} + b\bgrb\gamma^{\fieldChar^{i\cdot n}}+b^{2}  \nonumber    \\
&& \tab \bglb ~~ \textrm{obtained by substituting the expansion of ~}  \gamma^{(1+\fieldChar^{i\cdot n})} \textrm{~ from (\ref{Non-square-product-of-normal-basis-elements})}  ~~\bgrb \nonumber   \\
&& = ~~  \bglb (i\cdot h)^{-1}-b \bgrb \bglb \delta^{-1}+b \bgrb- \bglb (i\cdot h)^{-1}+b \bgrb \bglb \delta^{-\fieldChar^{i\cdot n}}+b \bgrb + b^{2} \nonumber  \\
&& = ~~  \bglb (i\cdot h)^{-1}-b\bgrb \delta^{-1} - \bglb (i\cdot h)^{-1}+b\bgrb \delta^{-\fieldChar^{i\cdot n}} - b^{2}
\label{Non-square-product-of-normal-basis-elements-with-redundant-representation}
\end{eqnarray}
Now, observing that 
$\sum_{i = 0}^{\fieldChar-1} \delta^{-\fieldChar^{i\cdot n}}$
$ = $
$\sum_{i = 0}^{\fieldChar-1} (\gamma-b)^{\fieldChar^{i\cdot n}}$
$ = $
$\sum_{i = 0}^{\fieldChar-1} \gamma^{\fieldChar^{i\cdot n}}$
$ = $
$-\alpha^{-1}$,
the following elaborations hold, for $\delta^{-2}$:
\begin{eqnarray*}
&& \ltab \ltab \ltab \ltab \ltab -\alpha^{-1} \delta^{-1}  ~~ = ~~
 \delta^{-1} \sum_{i = 0}^{\fieldChar-1} \delta^{-\fieldChar^{i\cdot n}} ~~ = ~~ \delta^{-2} +  \sum_{i = 1}^{\fieldChar-1} \delta^{-(1+\fieldChar^{i\cdot n})}    \\
&& \ltab \ltab = ~~  \delta^{-2} +  \sum_{i = 1}^{\fieldChar-1}\left [ ~ \bglb (i\cdot h)^{-1}-b \bgrb \delta^{-1} - \bglb (i\cdot h)^{-1} + b\bgrb \delta^{-\fieldChar^{i\cdot n}} - b^{2} ~ \right]    \\
&& \ltab \ltab   = ~~  \delta^{-2}  + b \delta^{-1} + b\bglb \alpha^{-1}+\delta^{-1}\bgrb + \sum_{i = 1}^{\fieldChar-1}\left [ (i\cdot h)^{-1}  \delta^{-1} -(i\cdot h)^{-1} \delta^{-\fieldChar^{i\cdot n}} - b^{2} ~ \right]    \\
&& \ltab \ltab   = ~~  \delta^{-2}  + b \delta^{-1} + b\bglb \alpha^{-1}+\delta^{-1}\bgrb + \sum_{i = 1}^{\fieldChar-1}\left [ (i\cdot h)^{-1}  \delta^{-1} -(i\cdot h)^{-1} \delta^{-\fieldChar^{i\cdot n}}  ~ \right]  + b^{2}\,,
\end{eqnarray*}
where the simplifications  $\sum_{i =1}^{\fieldChar-1}1  = -1$ and
 $ \sum_{i = 1}^{\fieldChar-1}\delta^{-\fieldChar^{i\cdot n}} =
-\alpha^{-1}-\delta^{-1} $ are performed in the last couple of expressions.
Thus, 
\begin{eqnarray*}
&&\ltab \ltab \ltab \ltab \delta^{-2} ~~ =  ~~  - ~ b\bglb b + \alpha^{-1}\bgrb ~ - ~ \left [ 2b + \alpha^{-1} + \sum_{i = 1}^{\fieldChar-1} (i\cdot h)^{-1} \right ] \delta^{-1} ~ + ~   \sum_{i = 1}^{\fieldChar-1} (i\cdot h)^{-1} \delta^{-\fieldChar^{i\cdot n}}  \\
&& \ltab \ltab = \left\{
\begin{array}{lcl}
- ~ b\bglb b + \alpha^{-1}\bgrb ~ - ~ \bglb 2b+\alpha^{-1}+1)\delta^{-1} ~  + ~ \sum_{i = 1}^{\fieldChar-1} (i\cdot h)^{-1} \delta^{-\fieldChar^{i\cdot n}}  & , & \textrm{for~} \fieldChar = 2, \textrm{~and}\\
~ - ~ b\bglb b +  \alpha^{-1}\bgrb ~ - ~ \bglb 2b+\alpha^{-1})\delta^{-1}  ~ + ~ \sum_{i = 1}^{\fieldChar-1} (i\cdot h)^{-1} \delta^{-\fieldChar^{i\cdot n}} & , & \textrm{for~} \fieldChar \geq 3 
\end{array} \right .
\end{eqnarray*}
If $\fieldChar \geq 3$, then $\sum_{i = 1}^{\fieldChar-1} (i\cdot h)^{-1} = 0$.
Looking at only the expression for $\delta^{-2}$, it may be tempting to set
$b = -\alpha^{-1}$,  but the expression for $\delta^{-(1+\fieldChar^{i\cdot n})}$,
elaborated in (\ref{Non-square-product-of-normal-basis-elements-with-redundant-representation}),
forbids this particular choice. Instead, it is preferable to let $b = 1$ and
include it in the multiplication table as an extra (redundant) table. Now,
in order to ensure that $\alpha^{-1}$ may have been constructed just as $\delta^{-1}$,
it must be checked that $\sum_{l = 0}^{n\fieldChar-1} \delta^{\fieldChar^{l}} \neq 0$,
where $n$ is the dimension of $\scalars$ over $\finiteIntegerField$, 
for $x^{\fieldChar}-x-\delta$ to be irreducible over $\scalars$. The minimal polynomial of
$\gamma$ over $\scalars$ is obtained from the equation
 $\gamma^{\fieldChar} +\alpha^{-1} \gamma^{\fieldChar-1} -\alpha^{-1} = 0$, and substituting
 $\gamma = \bglb \delta^{-1}+b \bgrb$ and multiplying both sides by $\delta^{\fieldChar}$,
 it may be found that
  $\bglb 1+b\delta\bgrb^{\fieldChar} +\alpha^{-1}\delta  \bglb 1+b\delta\bgrb^{\fieldChar-1} -\alpha^{-1}\delta^{\fieldChar} = 0$.
   The coefficient of $\delta^{\fieldChar}$ is $\bglb b^{\fieldChar}+\alpha^{-1}b^{(\fieldChar-1)}-\alpha^{-1} \bgrb$, 
  which cannot vanish for any $b \in \scalars$, by virtue of irreducibility of $x^{\fieldChar}-x-\alpha$ over $\scalars$,
  and the coefficient of $\delta^{(\fieldChar-1)}$ is $ \bglb(\fieldChar-1) b^{(\fieldChar-2)}\alpha^{-1}
  \bgrb \neq 0$,  for $b \neq 0$, and hence
 \begin{eqnarray*}
&& \shiftleft  \sum_{i = 0}^{\fieldChar-1} \delta^{\fieldChar^{i\cdot n}}
   ~~  = ~~  - \bglb b^{\fieldChar}+\alpha^{-1}b^{(\fieldChar-1)}-\alpha^{-1} \bgrb^{-1}
   \bglb(\fieldChar-1) b^{(\fieldChar-2)}\alpha^{-1}  \bgrb    \\
&& = ~~    
   b^{(\fieldChar-2)} \bglb ~~ \alpha \cdot \bglb b^{\fieldChar}+\alpha^{-1}b^{(\fieldChar-1)}-\alpha^{-1} \bgrb ~~ \bgrb^{-1} 
   ~~  = ~~  b^{(\fieldChar-2)} \bglb b^{\fieldChar}\alpha + b^{(\fieldChar-1)} - 1 \bgrb^{-1} \\
&& = ~~       b^{-2}  \alpha^{-1} ~~ ,\tab \tab \textrm{whenver ~} b \in \finiteIntegerField \backslash \{ 0\} 
\end{eqnarray*} 
  Now, $\delta^{-1}$ is a normal element in $\extensionField$ as an extension field of $\scalars$
  and the polynomial $x^{\fieldChar}-x-\delta$ becomes irreducible over $\extensionField$,
  as an extension field of $\scalars$, by Theorem \ref{Normal-Basis-of-Sec-Artin-Schreier-Extension-of-Finite-Fields-2},
  whenever $b \in \finiteIntegerField \backslash \{0\}$.
  A simpler proof could be as follows: 
  \begin{eqnarray*}  
&&  \ltab \ltab  \sum_{l = 0}^{n\fieldChar-1} \delta^{\fieldChar^{l}}  ~~ =  ~~
   \sum_{j = 0}^{n-1} \sum_{i = 0}^{\fieldChar-1} \delta^{\fieldChar^{i\cdot n + j}} 
    ~~  = ~~ \sum_{j = 0}^{n-1} \sum_{i = 0}^{\fieldChar-1} \delta^{\fieldChar^{i\cdot n } \cdot \fieldChar^{j}} 
   ~~ = ~~\sum_{j = 0}^{n-1} \left[ \sum_{i = 0}^{\fieldChar-1} \delta^{\fieldChar^{i\cdot n} } \right]^{\fieldChar^{j}} \\ 
&&  = ~~  \sum_{j = 0}^{n-1} \bglb  b^{-2}  \alpha^{-1} \bgrb^{\fieldChar^{j}} ~~ ,\tab \tab \textrm{whenver ~} b \in \finiteIntegerField \backslash \{ 0\}  \\
\end{eqnarray*}
and the last expression is the trace of $\bglb  b^{-2}\alpha^{-1} \bgrb$,
which cannot vanish, if $\alpha^{-1}$ is a normal element in $\scalars$,
as an extension field of $\finiteIntegerField$, for appropriate choice $b 
\in \finiteIntegerField \backslash \{0\}$.
In particular, if $b = 1$, the conjunct hypothesis that the trace of $\alpha$ is nonzero, together with
that $\alpha^{-1}$ is a normal element in $\scalars$ as an extension field of $\finiteIntegerField$,
are stated as part of the induction hypothesis.  It is still required to ensure that $\delta^{-1}$
is a normal element in $\extensionField$ over $\finiteIntegerField$, in order to ensure that the
induction hypothesis carries over to the nest step, making it necessary to apply Part 2 of
Theorem \ref{Normal-Basis-of-Sec-Artin-Schreier-Extension-of-Finite-Fields-1}.
 Insofar as the tedious computational details of 
 Theorem \ref{Normal-Basis-of-Sec-Artin-Schreier-Extension-of-Finite-Fields-2}
 are concerned, the simplifications discussed can be applied.
  Thus,   the following induction can be applied, for multiple Artin-Schreier extensions:
  initially, the polynomial $x^{\fieldChar}-x-1$ is irreducible over $\finiteIntegerField$,
   and $\alpha$ can be taken to be $1$,
  because this choice allows the element $1$ to be included in the initial basis,
  and as induction step, for a given normal element $\alpha^{-1}$, with the polynomial
  $x^{\fieldChar}-x-\alpha$  being irreducible over $\scalars$, a normal element
  $\delta^{-1} \in \extensionField$,  with the additional requirement that
  $x^{\fieldChar}-x-\delta$ is irreducible over $\extensionField$, is constructed.
  The inclusion of $1$ in the initial basis allows the multiplication tables
  to be represented succinctly, for efficient computation. The representation 
  of the product is not necessarily only in terms of linear combinations of the
  normal basis elements, making room for a redundant representation, with $1$ as
  the extra (redundant) element, being utilized for such a purpose.
  The computational speedup and succinctness of the representation
  cannot be overstated.
  
\textcolor{secheader}{
\begin{corollary}
 \label{Corollary-to-Normal-Basis-of-Sec-Artin-Schreier-Extension-of-Finite-Fields-1-and-2}
  {\bf {\textsf{~~(Inductive Construction of Near Optimal Normal Bases of Multiple Artin-Schreier-Extensions)~~}}} 
    Let $\scalars$ be a finite dimensional extension field of $\finiteIntegerField$, of vector space dimension $nn$, for some positive integer $n$. Let $\alpha \in \scalars \backslash  \{0\}$  be such that the polynomial $ x^{\fieldChar}-x- \alpha$ is an irreducible Artin-Schreier polynomial over $\scalars$, and $\alpha^{-1}$  is a normal element in $\scalars$ over
   the subfiled of $\fieldChar^{m}$  elements, and also a normal element in $\scalars$ over the prime field $\finiteIntegerField$, and let  $\extensionField = \scalars[\beta]/\bglb\beta^{\fieldChar}-\beta-\alpha\bgrb$.  Then, the element $\delta^{-1} = \beta^{-1}-1$  is a normal element in $\extensionField$ as an extension field of $\scalars$,  and also a normal element in 
   $\extensionField$ as an extension field of the prime field $\finiteIntegerField$,  and the polynomial  $ x^{\fieldChar}-x-\delta$ is an irreducible Artin-Schreier polynomial. Moreover, the multiplication table of $\delta^{-1}$ is nearly optimal, with the inclusion of $1$ as an extra (redundant) element, if so needed.
\end{corollary}
}
\proof. Follows form the discussion in the preceding paragraph, with $b = 1$.   			\qed

\textcolor{titlecolor}{
\section{Normal Bases and Their Optimality in Other Finite Field Extensions}
}
Let $\fieldChar$ and $\primeInt$ be distinct prime numbers, 
$l$ be a positive integer, such that $\primeInt\, \vert \, (\fieldChar^{l}-1)$,
and let $r$ be the largest positive integer such that 
$\primeInt^{r} \, \vert \, (\fieldChar^{l}-1)$.
Let $\scalars$ be a finite field containing $\fieldChar^{l}$ elements, and
$\xi \in \scalars$ be a primitive $\primeInt^{r}$-th root of $1$. Then,
for every positive integer $s$,  the polynomial
$x^{\primeInt^{s}} -\xi \in \singlevariablepolynomials{\scalars}{x}$ is irreducible,
as can be inferred from the tower of extensions
$\basel{x^{\primeInt}}{1} = \xi$ and $\basel{x^{\primeInt}}{i} = \basel{x^{\primeInt}}{i-1}$,
for $2 \leq i \leq s$, whenever $s \geq 2$ \cite{Shoup1990, Shoup1994}.
In the sequel, let $s$ be a positive integer, such that $s \leq r$.
Let $\extensionField$
$ = $
$\singlevariablepolynomials{\scalars}{\alpha}/(\alpha^{\primeInt^{s}}-\xi)$,
and $\fieldChar^{l}-1 = m \primeInt^{r}$, with $\gcd(m,\, \primeInt) = 1$.
Now, $\alpha^{\fieldChar^{l}-1} = \xi^{m\primeInt^{r-s}} = \zeta$,
for some $\zeta \in \scalars$, such that $\zeta \neq 1$ and
$\zeta^{\primeInt^{s}} = 1$. If $\primeInt \, \vert \, (\fieldChar-1)$, then 
$\zeta \in \basel{\Integers}{\fieldChar} \backslash \{1\}$.
The least positive integer $k$, such that $\zeta^{k} = 1$,
is $k = \primeInt^{s}$. Thus,
$\alpha^{\fieldChar^{l}} - \zeta \alpha = 0$,
and hence, 
$\alpha^{\fieldChar^{i\cdot l}} - \zeta^{i} \alpha = 0$, for $1 \leq i \leq \primeInt^{s}-1$.
Let $\beta = \alpha-b$, for some $b \in \scalars \backslash \{0\}$, such that
$b^{\primeInt^{s}} \neq \xi$. Then,
$\beta^{\fieldChar^{i\cdot l}}$
$ = $
$\alpha^{\fieldChar^{i\cdot l}} - b$
$ = $
$\zeta^{i} \alpha  - b$
$ = $ 
$\zeta^{i} (\beta+b) - b$,
and hence,
$\beta^{\fieldChar^{i\cdot l}}$
$ - $
$\zeta^{i} \beta$
$ = $
$(\zeta^{i}-1)b$,
where  $\zeta^{i},\, b \in \scalars \backslash \{0\}$,
$\zeta^{i} \neq 1$,
for $1 \leq i \leq \primeInt^{s}-1$.
Now, the element $\gamma = \beta^{-1}$ satisfies
$(\zeta^{i}-1)b \gamma^{1+\fieldChar^{i\cdot l}}$
$ = $
$\gamma$
$ - $
$\zeta^{i} \gamma^{\fieldChar^{i\cdot l}}$,
with $\zeta^{i},\, b, \, (\zeta^{i}-1)b \in \scalars \backslash \{0\}$,
for $1 \leq i \leq \primeInt^{s}-1$.
It is convenient and preferable to choose $b = 1$.
The expansion for the term $\gamma^{2}$
must be found from the formula 
$\bgls\sum_{i = 0}^{\primeInt^{s}-1} \gamma^{\fieldChar^{i\cdot l}} \bgrs \cdot \gamma$
$ = $
$\gamma^{2} + \sum_{i = 1}^{\primeInt^{s}-1} \gamma^{1+\fieldChar^{i\cdot l}} $ 
$ = $
$\gamma^{2} + \sum_{i = 1}^{\primeInt^{s}-1}(\zeta^{i}-1)^{-1} b^{-1} \gamma$
$ - \sum_{i = 1}^{\primeInt^{s}-1}(\zeta^{i}-1)^{-1} b^{-1} \zeta^{i}\gamma^{\fieldChar^{i\cdot l}} $.
The minimal polynomial of $\alpha$ in
$\singlevariablepolynomials{\scalars}{x}$ is $x^{\primeInt^{s}}-\xi$,
that of $\beta$ is $(x+b)^{\primeInt^{s}}-\xi$ and that of $\gamma$ is 
$(1+bx)^{\primeInt^{s}}-\xi x^{\primeInt^{s}}$,
except for multiplication by $\bglb b^{\primeInt^{s}}-\xi\bgrb$.
From the minimal polynomial, it follows that
$\sum_{i = 0}^{\primeInt^{s}-1} \gamma^{\fieldChar^{i\cdot l}}$
$ = $
$-\primeInt^{s} b^{\primeInt^{s}-1}\bglb b^{\primeInt^{s}}-\xi\bgrb^{-1}$.
Now, in order to show that $\gamma$ becomes a normal element in $\extensionField$,
it may be recalled that $\beta^{\fieldChar^{i\cdot l}}$
$ = $
$\alpha^{\fieldChar^{i\cdot l}} - b$
$ = $
$\zeta^{i} \alpha  - b$, and hence
$\gamma^{\fieldChar^{i\cdot l}}$
$ = $
$(\zeta^{i} \alpha  - b)^{-1}$
$ = $
$\alpha^{-1}(\zeta^{i} - \alpha^{-1} b)^{-1}$,
for $0 \leq i \leq \primeInt^{s}-1$.
The identity $ 1 - x^{\primeInt^{s}} $
$ = $
$\zeta^{i\cdot \primeInt^{s}} -x^{\primeInt^{s}}$
$ = $ 
$(\zeta^{i}-x) \sum_{j = 0}^{\primeInt^{s}-1} \zeta^{i\cdot (\primeInt^{s}-j-1)}x^{j}$
holds, with $x$ as an indeterminate. Substituting $\tau$ for $x$,
where $\tau$ is an element in the algebraic closure
of $\scalars$, such that $\tau^{\primeInt^{s}} \neq 1$,
it follows that
$(\zeta^{i}-\tau)^{-1}$
$ = $
$(1 - \tau^{\primeInt^{s}})^{-1}  \sum_{j = 0}^{\primeInt^{s}-1} \zeta^{i\cdot (\primeInt^{s}-j-1)}\tau^{j}$,
and, with $\tau = \alpha^{-1}b$, 
it can be readily checked that
$\tau^{\primeInt^{s}}$
$ = $
$ \bglb \alpha^{-1}b \bgrb^{\primeInt^{s}}$
$ = $
$ \alpha^{-\primeInt^{s}} \cdot b^{\primeInt^{s}} $
$ = $
$ \xi^{-1} \cdot  b^{\primeInt^{s}}  \neq 1$,
and hence 
$\gamma^{\fieldChar^{i\cdot l}}$
$ = $
$\alpha^{-1}(\zeta^{i} - \alpha^{-1} b)^{-1}$
$ = $
$\alpha^{-1}(1 - (\alpha^{-1} b)^{\primeInt^{s}})^{-1}  \sum_{j = 0}^{\primeInt^{s}-1} \zeta^{i\cdot (\primeInt^{s}-j-1)}(\alpha^{-1} b)^{j}$, for $0 \leq i \leq \primeInt^{s}-1$.
Now, for $0 \leq i, \, j, \, k \leq \primeInt^{s}-1$, if $j \neq k$,
then $\zeta^{(\primeInt^{s}-j-1)} \neq \zeta^{(\primeInt^{s}-k-1)}$,
and hence, the $\primeInt^{s} \times \primeInt^{s}$ matrix with
$\zeta^{i\cdot (\primeInt^{s}-j-1)}$ as the entry in the $(i+1)$-th row and $(j+1)$-th
column is invertible. Thus, the linear span of the set with coefficients in $\scalars$
of the set $\{\gamma^{\fieldChar^{i\cdot l}}\, : \, 0 \leq i \leq \primeInt^{s}-1\}$
is the same as that of the set 
$\{ \alpha^{-j} \, : \, 0 \leq j \leq \primeInt^{s}-1 \}$,
and $\gamma$ is a normal basis element, with a nearly optimal multiplication table.
If $\zeta,\, b \in \basel{\Integers}{\fieldChar} \backslash \{0\}$,
which, in turn, requires that $\primeInt\, \vert \, (\fieldChar-1)$,
then most of the coefficients in the multiplication table are elements
in the prime field $\basel{\Integers}{\fieldChar}$.

\textcolor{titlecolor}{ 
\section{Summary and Conclusions}
}

In this paper, the minimal polynomial of a generating element of a finite dimensional extension over a finite field is described.  From the generating element of an extension, normal basis elements can be easily obtained. The minimal polynomials of such normal basis elements can be found from that of the generating element. It is then shown that the multiplication tables satisfy a near optimality property.

\textcolor{titlecolor} {
 \section*{Acknowledgements}
 }
 \textcolor{secheader}{
 The authors are greatly indebted to Professor Michael Oser Rabin, taking the cue from Professor Amir Pnueli.
The authors have also had fruitful discussions with Professor Gary Mullen, Professor Neal Koblitz and Professor
Igor Shparlinski, during the preparation of the work.
}

 \textcolor{refscolor}{
}
\end{document}